\RequirePackage{fix-cm}
\documentclass[twocolumn]{svjour3}          
\smartqed  
\usepackage{graphicx}
\usepackage{xfrac}
\usepackage[utf8]{inputenc}
\usepackage{cancel}
\usepackage{mathtools,cases}
\usepackage{amsmath}
\usepackage{amsbsy}
\usepackage{csquotes}
\usepackage[space]{cite}
\usepackage{amsfonts}
\usepackage{url}
\usepackage[dvipsnames]{xcolor}
\usepackage{booktabs}
\usepackage{dblfloatfix} 
\usepackage[overload]{empheq}
\usepackage{subcaption}
\DeclareCaptionSubType*[Alph]{table}
\DeclareCaptionLabelFormat{mystyle}{Table~\bothIfFirst{#1}{ }#2}
\newcommand{\norm}[2]{ \Vert #1 \Vert_{#2}}
\begin{document}

\title{A Singularity Removal Method for Coupled 1D-3D Flow Models}

\author{Ingeborg G. Gjerde         \and
        Kundan Kumar \and 
        Jan M. Nordbotten
}

\institute{I. G. Gjerde \at
              Department of Mathematics, University of Bergen \\
              \email{ingeborg.gjerde@uib.no}           
           \and
           K. Kumar \at
	Department of Mathematics, University of Karlstad \\
	          \email{kundan.kumar@karlstad.se}         \\   %
	           Department of Mathematics, University of Bergen \\
	                     \and
           J. M. Nordbotten \at
	Department of Mathematics, University of Bergen \\
	          \email{jan.nordbotten@uib.no}           
}

\date{Received: date / Accepted: date}

\maketitle

\begin{abstract}
In reservoir simulations, the radius of a well is inevitably going to be small compared to the horizontal length scale of the reservoir. For this reason, wells are typically modelled as lower-dimensional sources. In this work, we consider a coupled 1D-3D flow model, in which the well is modelled as a line source in the reservoir domain and endowed with its own 1D flow equation. The flow between well and reservoir can then be modelled in a fully coupled manner by applying a linear filtration law. 

The line source induces a logarithmic type singularity in the reservoir pressure that is difficult to resolve numerically. We present here a singularity removal method for the model equations, resulting in a reformulated coupled 1D-3D flow model in which all variables are smooth. The singularity removal is based on a solution splitting of the reservoir pressure, where it is decomposed into two terms: an explicitly given, lower regularity term capturing the solution singularity and some smooth background pressure. The singularities can then be removed from the system by subtracting them from the governing equations. Finally, the coupled 1D-3D flow equations can be reformulated so they are given in terms of the well pressure and the background reservoir pressure. As these variables are both smooth (i.e. non-singular), the reformulated model has the advantage that it can be approximated using any standard numerical method. The reformulation itself resembles a Peaceman well correction performed at the continuous level.

\end{abstract}

\section{Introduction}
\label{intro}
Accurate well models are of critical importance for reservoir simulations. The well constitutes the driving force for reservoir flow, in addition to being the main access point of information about its state. The major challenge of well modelling is that of scale disparity; a well has a radius of $\sim 10$ cm, while the reservoir might extend several kilometres in the horizontal plane. From a computational viewpoint, this makes it exceedingly expensive to resolve the well as a 3D object in the grid representing the reservoir. For this reason, wells are typically modelled using either zero-dimensional (0D) point sources or (1D) line sources.

\begin{figure}[t]
    \begin{center}
\includegraphics[width=0.2\textwidth]{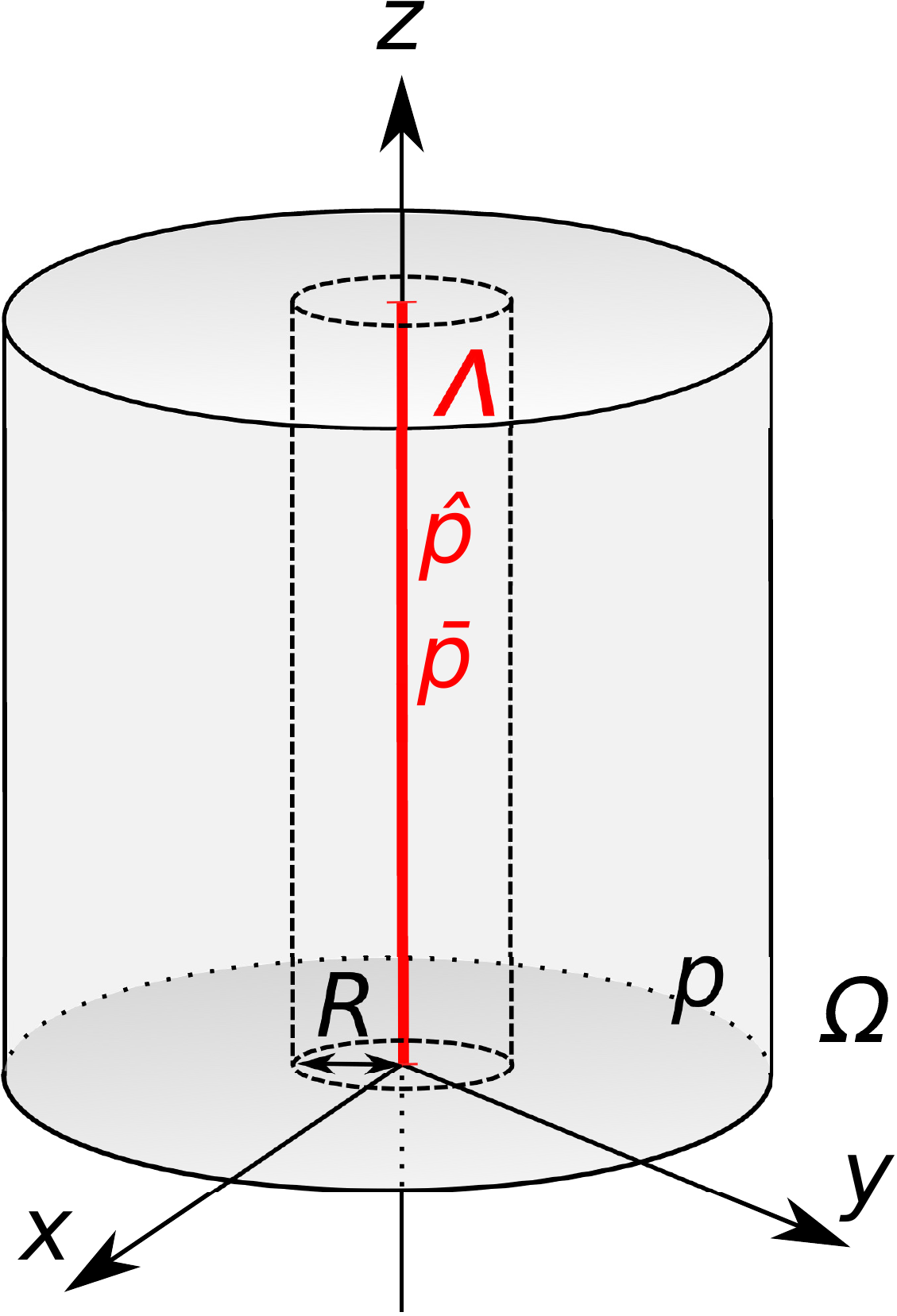}
     \end{center}
     \caption{A 1D domain $\Lambda$ embedded in a 3D domain $\Omega$ representing the reservoir. The reservoir domain $\Omega$ is allowed to be arbitrarily shaped. The well is considered to be a thin cylinder of radius $R \ll \text{size}(\Omega)$. For this reason, the radial and angular components of the well pressure $\hat{p}$ are ignored, so that it can be described as a 1D variable $\hat{p} = \hat{p}(s)$.}
   \label{fig:illustration}
\end{figure}

In this work, we take as a starting point the coupled 1D-3D flow model
\begin{subequations}
\begin{align}
 \mathbf{q} + \frac{\kappa}{\mu}  \nabla p =& \,  0 & \text{in} \, \Omega,\label{eq:basic3D1} \\
 \nabla \cdot \mathbf{q} =&  \, \beta \left( \hat{p}- \bar{p}\right) \delta_{\Lambda}  & \text{in} \, \Omega,\label{eq:basic3D2}  \\
 \hat{q} + \frac{\hat{\kappa}}{\mu} \frac{\mathrm{d} \hat{p}}{\mathrm{d} s} =& \, 0 & \text{in } \, \Lambda,\label{eq:basic1D1}\\
\frac{\mathrm{d} \hat{q}}{\mathrm{d} s}   =& - \hat{\beta} \left( \hat{p}- \bar{p}\right)  & \text{in} \, \Lambda ,\label{eq:basic1D2} 
\end{align}
\end{subequations}
where $\Omega \subset \mathbb{R}^3$ denotes the reservoir domain and $\Lambda = \cup_{w=1}^{\text{wells}} \Lambda_w \subset \mathbb{R}^1$ a collection of line segments each representing a well. The 1D domain is parametrized by its arc-length $s$. The parameters $\kappa$, $\hat{\kappa}$ and $\mu$ denote reservoir permeability, well permeability and fluid viscosity, respectively, and are assumed to be positive and constant. The variables $p$ and $\mathbf{q}$ denote fluid pressure and flux in the reservoir, $\hat{p}$ and $\hat{q}$ fluid pressure and flux in the well, and $\bar{p}$ the reservoir pressure averaged over the surface of the borehole $r=R$:
\begin{align}
\bar{p}(z, R) = \frac{1}{2\pi R} \int_0^{2\pi} p(R, z, \theta) \mathrm{d} \theta,\label{eq:avg}
\end{align}
as is illustrated in Figure \ref{fig:illustration}.

Physically, equations \eqref{eq:basic3D1}-\eqref{eq:basic3D2} describe a Darcy-type flow in the reservoir domain $\Omega$, and equations \eqref{eq:basic1D1}-\eqref{eq:basic1D2} a Poiseuille-type flow in the well. The latter is a 1D flow equation, where the radial and angular components have been neglected. For a description of this model reduction method for the well flow, we refer to the work of Cerroni et al. in \cite{Cerroni2019}. The mass flux $q$ between reservoir and well is modelled using a linear filtration law,
\begin{align}
q = \beta (\hat{p}- \bar{p}),
\end{align}
which states that the connection flow between them is proportional to their pressure difference. The proportionality coefficients $\beta, \hat{\beta} \in C^{1}(\Lambda_w)$ are assumed piecewise continuous and allowed to vary along the well. The wells are considered as concentrated line sources $\delta_\Lambda$ in the reservoir equation \eqref{eq:basic3D1}-\eqref{eq:basic3D2}, with the line sources defined in the following manner:
\begin{align}
 \int_\Omega f \delta_\Lambda \, \phi \, \mathrm{d} \Omega  &= \sum_{w=1}^{\text{wells}} \int_{\Lambda_w} f(s_w) \phi(s_w) \mathrm{d}s_w\label{eq:dirac}
\end{align}
for all $ \phi \in C^0(\Omega)$, with $s_w$ denoting the arc-length of line segment $\Lambda_w$.

Elliptic equations with line sources of the type \eqref{eq:dirac} have been used in a variety of applications, e.g., the modelling of 1D steel components in concrete structures \cite{llau2016} or the interference of metallic pipelines and bore-casings in electromagnetic modelling of reservoirs \cite{weiss2017}. A coupled 1D-3D heat transfer problem was considered in the context of geothermal energy in \cite{alkhoury2005}, where it was used to model heat exchange between (3D) soil and a (1D) pipe. Coupled 1D-3D flow models have also been studied in the context of biological applications, such as the efficiency of cancer treatment by hyperthermia \cite{nabil2016}, the efficiency of drug delivery through microcirculation \cite{Cattaneo2014, zunino2018}, and the study of blood flow in the vascularized tissue of the brain \cite{Reichold2009, Grinberg2011}. In this work, we restrict ourselves to considering its application in the context of reservoir modelling. 

The main challenge with the coupled 1D-3D flow problem is that the line source induces the reservoir pressure to be singular, thereby making its analysis and approximation non-standard. Typically, reservoir simulations are performed using finite volume methods. The discretized form of the coupling in \eqref{eq:basic3D1}-\eqref{eq:basic1D2} is then given by
\begin{align}
q = \beta (\hat{p}-p_K),
\end{align}
where $p_K$ denotes the average pressure in the grid block containing the well. Due to the singularity, $p_K$ will not be representative of the reservoir pressure at the bore-hole; this is typically accounted for by multiplying $\beta$ with a well index $J$. A correction of this type was first developed by Peaceman in \cite{peaceman1978}, where he considered the two-point flux approximation method on uniform, square grids when the well is aligned with one of its axes. Via an analytic solution valid for simplified cases, he gave a well index depending on the equivalent radius of the well, i.e., the radius at which the reservoir pressure equals the well block pressure. The equivalent radius depends, among other factors, on the discretization scheme, placement of the well relative to the mesh, and reservoir permeability. The problem of finding appropriate well indexes has been treated in a multitude of works; Peaceman himself treated an extension of his method to non-square grid-blocks and anisotropic permeability \cite{Peaceman-rectangle}. The extension to more generalized grids was treated by e.g. Aaavatsmark in \cite{Aavatsmark2016-0, Aavatsmark2016,  Aavatsmark2016-2}, to more generalized flow models by e.g. Ewing in \cite{Ewing1999}, and to more generalized discretization schemes by e.g. Chen et al. in \cite{Chen2009}. Many authors have contributed to the extension to generalized well placements, we mention here the work of King et. al in \cite{king}, Aavatsmark in \cite{aavatsmark2003index}, and of special relevance to our work, that of Wolfsteiner et al. in \cite{Wolfsteiner2003} and Babu et al. in \cite{babu}. 

In this work, we take a different approach, in which the singularities are explicitly removed from the governing equations. We start by showing that the reservoir pressure $p$ admits a splitting
\begin{align}
p = \sum_{w=1}^{\text{wells}} E \left(\beta(\hat{p}-\bar{p})  \right) \Psi_w G_w  + v,
\label{eq:decomp1}
\end{align}
where $G_w$ is a given logarithmic function that captures the near-well behaviour of the reservoir pressure, $E$ is an extension operator $E \colon H^2(\Lambda) \rightarrow H^2(\Omega)$, $\Psi_w$ some smooth cut-off function, and $v \in H^2(\Omega)$ some higher-regularity remainder term. The key point here is that the singular nature of the solution is explicitly captured by the logarithmic terms $G_w$. With the splitting \eqref{eq:decomp1} in hand, we can therefore remove the singular terms from the system by straightforward subtraction. Finally, we reformulated coupled 1D-3D flow model can then be reformulated so it is given with respect to the  high-regularity variables $\hat{p}$ and $v$. The main contribution of this article is the reformulation of the coupled 1D-3D flow model into equations \eqref{eq:reform-3D}-\eqref{eq:reform-1D-bc}, for which the solution is smooth (non-singular). On a practical level, this means the solution can be approximated using any standard numerical method.

The technique of removing singularities is commonly known for point sources; we refer here to \cite[p. 14]{ewing-book} for a more in depth explanation. It has previously been studied in the context of reservoir models by e.g. Hales, who used it to improve  well modelling for 2D reservoir models \cite{hales}. A splitting of the type \eqref{eq:decomp1} was introduced by Ding in \cite{Ding2001} for the point source problem, where it was used to formulate grid refinement strategies. We are, to the best of our knowledge, the first to formulate a singularity removal method for the coupled 1D-3D flow problem. Central to this method is the construction of a function $G_w$ capturing the solution singularity; we use here a function $G_w$ found by integrating the Green's function for the reservoir equations \eqref{eq:basic3D1}-\eqref{eq:basic3D2} over the line $\Lambda$; we refer here to our earlier work in \cite[Section 3.2]{Gjerde2018}. This use of Green's functions to construct analytical and semi-analytical well models has a rich history.  Of special relevance to our work, we mention that of Wolfsteiner et al. and Babu et al. in \cite{babu, Wolfsteiner2003}, in which the Green's function was used to construct analytical solutions with which to calculate the well index $J$. More recently, Nordbotten et al. used Green's functions to construct analytical models to estimate leakage of CO$_2$ stored in geological formations \cite{Nordbotten2009}. 

The singularity removal, and subsequent reformulation of the model in terms of the smooth variables $v$ and $\hat{p}$, is similar to the Peaceman well correction in that it leads to an alteration of the inflow parameter $\beta$. We discuss this in more detail in Section \ref{sec:peaceman}. It differs, however, in that it works on the continuous level. It is therefore easily adapted to different discretization methods, generalized well placements within the domain and different types of boundary conditions. Moreover, since our method gives an explicit representation of the logarithmic nature of the solution, it allows us to accurately represent the reservoir pressure in the whole domain (including in the near-vicinity of the well).

In our presentation of the method, we limit ourselves to considering a linear reservoir equation with constant, scalar-valued permeabilities and Poiseuille flow in the well. The latter restriction is not critical to the methodology; the well equation could for example be taken non-linear as long as the well pressure remains sufficiently regular. To be more precise, the method requires $\hat{p}$ to be piecewise $C^1$ on $\Lambda$. As for the reservoir equation, the reservoir pressure could be replaced with a potential expression $\phi$ so that the effect of gravity can be included. The singularity removal and reformulation can be extended to handle spatially varying, scalar-valued permeabilities as shown in \cite{Gjerde2018}. For an extension to tensor-valued permeabilities and non-linear reservoir equations, we suggest using the solution splitting in \eqref{eq:decomp1} to formulate a multiscale finite volume method such as in \cite{hamdi2006}, or a generalized finite element method \cite{gfem}, where the analytic functions capturing the solution singularity are used to enrich the set of basis functions.

For the discretization and numerical experiments, we consider herein the Galerkin Finite Element (FE) method. The FE approximation of the line source problem was studied by D'Angelo in \cite{dangelo2012} by means of weighted Sobolev spaces, using similar techniques as those known for e.g. corner-point problems \cite{babuska1972}. D'Angelo proved that the approximation of the coupled 1D-3D flow problem \eqref{eq:basic3D1}-\eqref{eq:basic3D2} converges sub-optimally unless the mesh is sufficiently refined around the well. The sub-optimal convergence rates were found to be local to the line source by Köppl et al. in \cite{koppl2015}, meaning that they only pollute the pressure approximation inside the well block. However, this means the approximation of the coupled 1D-3D flow problem will suffer until the mesh parameter $h$ is smaller than the well radius $R$. In practice, one therefore needs a very fine mesh around the well for the FE approximation of \eqref{eq:basic3D1}-\eqref{eq:basic1D2} to converge. This makes the problem computationally expensive to solve. Different strategies have been proposed to remedy this, e.g., Kuchta et al. studied suitable preconditioners in \cite{miro2016-3D1D}. Holter et al. then applied this preconditioner to simulate flow through the microcirculature found in a mouse brain \cite{Holter2018}. An alternative coupling scheme was introduced  by Köppl et al. in \cite{koppl2016}, where the source term was taken to live on the boundary of the inclusions. The result is a 1D-(2D)-3D method where the approximation properties have been improved, at the expense of having to resolve the 2D boundary of the well.

The article is structured as follows. We start in Section  \ref{sec:notation} by defining the relevant function spaces for the problem. In Section \ref{sec:model}, we introduce in more detail the coupled 1D-3D flow model we take as a starting point. In Section \ref{sec:decomposition}, we show that the reservoir pressure $p$ admits a splitting into lower-regularity terms that capture the solution singularities, and a higher-regularity remainder term $v$. 
With the splitting in hand, the singularities can then be subtracted from the governing equations. The result is the reformulated coupled 1D-3D flow model \eqref{eq:reform-3D}-\eqref{eq:reform-1D-bc}, posed in terms of the smooth variables $\hat{p}$ and $v$. As the solutions then enjoy significantly improved regularity, this system can be approximated using standard numerical methods. The variational formulation and FE discretization of the reformulated problem are given in Sections \ref{sec:variational-form} and \ref{sec:disc}, respectively, and require only standard function spaces. In Section \ref{sec:peaceman}, we discuss how this discretization of the reformulated model resembles a Peaceman well correction. We then conclude the article with two numerical experiments, where we test the Galerkin FE method of both the standard and reformulated coupled 1D-3D flow model. We show that the singularity removal recovers optimal convergence rates on uniform meshes, i.e., without needing to refine the mesh around the well. Moreover, in a manner similar to altering the well index, it makes the approximation robust with respect to the ratio $R/h$. 

\section{Background and notation}
\label{sec:notation}
The purpose of this section is to introduce the appropriate function spaces for the coupled 1D-3D flow model. Let $H^{k}(\Omega)$ be the Sobolev space, 
\begin{align*}
H^{k}(\Omega) = \lbrace u \in L^2(\Omega) : D^\beta u \in L^2(\Omega) \text{ for } \vert \beta \vert \leq k \rbrace,
\end{align*}
with $\beta$ denoting a multi-index and $D^\beta$ the corresponding weak distributional derivative of $u$. $H^{k}(\Omega)$ is a Hilbert space endowed with inner product
\begin{align*}
(u, v)_{H^k(\Omega)} = \sum_{\vert \beta \vert \leq  k} \int_\Omega D^{\beta} u \, D^{\beta} v\, \mathrm{d} \Omega. 
\end{align*}
We use a subscript to denote the subspace of $H^k$ with zero trace on the boundary, $H^k_0$, i.e.,
\begin{align*}
H^k_0(\Omega) = \lbrace u \in H^k(\Omega) : u \vert_{\partial \Omega} = 0 \rbrace.
\end{align*} 

As we will see, the reservoir solution $p$ in \eqref{eq:basic3D1}-\eqref{eq:basic1D2} fails to belong to $H^1(\Omega)$ due to singular behaviour on $\Lambda$. For this reason, we consider also a weighted Sobolev space. To define it, let $ -1 < \alpha < 1$, and take $L^2_\alpha(\Omega)$ to denote the weighted Hilbert space consisting of measurable functions $u$ such that 
\begin{align*}
\int_\Omega u^2 r^{2\alpha} \mathrm{d} \Omega < \infty,
\end{align*}
where $r$ denotes the distance of a point to $\Lambda$, i.e., $r(\mathbf{x})=\text{dist}(\mathbf{x}, \Lambda)$. This space is equipped with the inner product 
\begin{align*}
(u, v)_{L_\alpha^2(\Omega)} = \int_\Omega r^{2\alpha} u v \, \mathrm{d} \Omega.
\end{align*}

For $\alpha>0$, the weight $r^\alpha$ has the power to dampen out singular behaviour in the function being normed; for $\alpha<0$, the weight function can induce or worsen already singular behaviour. We therefore have the relation $L^2_{-\alpha}(\Omega) \subset L^2(\Omega) \subset L^2_\alpha(\Omega)$ for $\alpha>0$. Letting now $H^{1}_\alpha(\Omega)$ be the Sobolev space
\begin{align*}
H^{1}_\alpha(\Omega) = \lbrace u \in L^2_\alpha(\Omega) : D^\beta u \in L^2_\alpha(\Omega) \text{ for } \vert \beta \vert \leq k \rbrace,
\end{align*}
we will later find that the reservoir pressure solving \eqref{eq:basic3D1}-\eqref{eq:basic1D2} belongs to $H^1_\alpha(\Omega)$ for $\alpha>0$.

A practical use of this space is found, for example, considering the logarithmic grading (refinement) that is often performed on a mesh around the well. The well introduces a logarithmic type singularity in the reservoir pressure that cannot be resolved using linear elements. Consequently, the convergence rate of standard numerical methods degrade using uniform meshes. Optimal convergence can be retrieved by a specific refinement of the mesh around the well \cite{apel2011, dangelo2012, Ding2001}. The exact convergence rates and mesh grading requirements are closely related to the weighted Sobolev space wherein the solution exists; in fact, the
graded mesh will be uniform with respect to the weight function $r^{\alpha}$. 

\section{Mathematical model}
\label{sec:model}
Here, we introduce in more detail the coupled 1D-3D equation we take as a starting point. Let $\Omega \subset \mathbb{R}^3$ denote a bounded domain describing a reservoir, with smooth boundary $\partial \Omega$. We consider here steady-state, incompressible Darcy flow
\begin{align}
\mathbf{q} = - \frac{\kappa}{\mu} \nabla p,
\end{align}
where $\mathbf{q}$ and $p$ denote reservoir flow and pressure, $\mu$ the fluid viscosity, and $\kappa$ a given positive and scalar permeability. We consider also a collection of wells, each considered to be a thin tube with fixed radius $R$ and centreline $\Lambda_w$. The centreline is parametrized by the arc length $s_w$. We denote by $\pmb{\tau}_{s_w}$ its normalized tangent vector. As the radius of the tube is small, we assume the radial and angular components of the well pressure can be neglected, meaning $\hat{p} \vert_{\Lambda_w}=\hat{p}(s_w)$. The well flow domain $\Lambda$ will then consist of a collection of line segments, $\Lambda = \cup_{w=1}^{\text{wells}} \Lambda_w$. We consider on this domain Poiseuille-type flow,
\begin{subequations}
\begin{align}
\hat{\mathbf{q}}_w &= \, - \frac{{R}^2}{8 \mu} \frac{\mathrm{d}\hat{p}}{\mathrm{d} s_w} \, \pmb{\tau}_{s_w}, \label{eq:poiseflow1}\\ 
\frac{\mathrm{d}\hat{\mathbf{q}}_w}{\mathrm{d} s_w}  &= -  \frac{q}{\pi R^2}, \label{eq:poiseflow2}
\end{align}
\end{subequations}
with $\hat{\mathbf{q}}_w$ and $\hat{p}_w$ denoting flow and pressure in the well and $q$ the linear mass flux into or out of the well. $\frac{\mathrm{d}}{\mathrm{d} s_w}$ denotes the derivative with respect to the tangent line, or equivalently, the projection of $\nabla$ along $\pmb{\tau}$, i.e., $\frac{\mathrm{d}}{\mathrm{d} s_w} = \nabla \cdot \pmb{\tau}_{s_w}$. As the fluid flux in the well has a fixed direction, it can be given as a scalar function $\hat{q}_w$, characterized by the property $\hat{\mathbf{q}}_w= \hat{q}_w \pmb{\tau}_w$. Note that the assumption of Poiseuille flow is not critical; \eqref{eq:poiseflow1} could for example contain certain non-linearities.

Letting now $\Lambda= \cup_{w=1}^\text{wells} \Lambda_w$ denote the collection of line segments $\Lambda_w$, the well pressure and flux can be written as 1D variables $\hat{p}, \hat{q} \colon \Lambda \rightarrow \mathbb{R}$. The well and reservoir flow can then be coupled together using a linear filtration law, which states that the mass flux $q$ between them is proportional to their pressure difference:
\begin{align}
q = 2\pi \lambda R f(\hat{p}, \bar{p}) \quad \text{where } f(\hat{p}, \bar{p})= \hat{p} - \bar{p}.
\end{align}
The mass flux is given as the rate of transfer per unit length, and the variable $\lambda \in C^2(\Lambda)$ denotes the permeability of the borehole lateral surface. It accounts for the fact that the well may not be in perfect contact with the reservoir, leading to a pressure drop across the borehole. Letting $\Delta p_{\text{skin}}$ denote this pressure drop, this can be expressed by the following relation: $q = 2\pi R \lambda \Delta p_{\text{skin}}$.

The pressure difference $f(\hat{p}, \bar{p})$ between well and reservoir uses an averaged value $\bar{p}(z; R)$ for the reservoir pressure given in \eqref{eq:avg}. This can be interpreted physically as the reservoir pressure averaged around the borehole. The flow in well and reservoir can be then modelled, in a fully coupled manner, by the set of equations
\begin{subequations}
\begin{align}
 \mathbf{q} + \frac{\kappa}{\mu}  \nabla p =& \,  0 & \text{in} \, \Omega,  \\
 \nabla \cdot \mathbf{q} =&  \, \beta f(\hat{p}, \bar{p}) \delta_{\Lambda}  & \text{in} \, \Omega,   \\
 p =& \, p_D & \text{on} \, \partial \Omega,   \\
 \hat{q} + \frac{\hat{\kappa}}{\mu} \frac{\mathrm{d} \hat{p}}{\mathrm{d} s}  =& \, 0 & \text{in} \, \Lambda, \\
\frac{\mathrm{d} }{\mathrm{d} s} \hat{q}   =& - \hat{\beta} f(\hat{p}, \bar{p})  & \text{in} \, \Lambda , \\
 \hat{p} =& \, \hat{p}_D & \text{on} \, \partial \Lambda,
\end{align}
\end{subequations}
where $\hat{\kappa} =  R^2/{8}$, $\beta=2\pi R \lambda$, $\hat{\beta} = \beta / \pi R^2$.  The functions $p_D \in C^2(\bar{\Omega})$ and $\hat{p}_D(\bar{\Lambda})$ denote given boundary data. The connection flow from well to reservoir is modelled by means of a generalized Dirac delta function $\delta_\Lambda$, which we understand in the sense of \eqref{eq:dirac}. Finally, this system can be reduced to its conformal form by eliminating the 1D and 3D fluxes:
\begin{subequations}
\begin{align}
\nabla \cdot \left( -\frac{\kappa}{\mu} \nabla p \right) =&  \beta f(\hat{p}, \bar{p}) \delta_{\Lambda} & \text{in} \, \Omega,\label{eq:3D-eq} \\
 p =& p_D & \text{on} \, \partial \Omega,\label{eq:3D-eq-bc} \\
 \frac{\mathrm{d} }{\mathrm{d} s} \left( - \frac{\hat{\kappa}}{\mu}\frac{\mathrm{d} }{\mathrm{d} s}  \hat{p} \right) =& - \hat{\beta} f(\hat{p}, \bar{p})   & \text{in} \, \Lambda,\label{eq:1D-eq}\\
 \hat{p} =& \hat{p}_D & \text{on} \, \partial \Lambda,\label{eq:1D-eq-bc}
\end{align}
\end{subequations}
with $f(\hat{p}, \bar{p})=\hat{p}-\bar{p}$.

\section{Splitting Properties of the Solution}
\label{sec:decomposition}
In this section, we will show that the line source in the right-hand side of \eqref{eq:3D-eq} introduces a particular structure to the solution of the coupled 1D-3D flow problem. We do this by means of a splitting technique, in which the reservoir pressure is split into a low regularity term that explicitly captures the singularity, and a regular component $v$ being the solution of a suitable elliptic equation. To start with, we discuss in detail the splitting when $\Lambda$ is assumed a single line segment aligned with the $z$-axis, $\frac{\kappa}{\mu}=1$ and the well outflow $q$ is a given function $f \in C^1_0(\Lambda)$. The splitting is then especially simple; this case therefore serves to illustrate the splitting method itself. We then generalize it in two steps, handling first an arbitrary line segment and $\frac{\kappa}{\mu}\neq 1$, and finally the coupling between reservoir and well. Finally, we use the splitting to reformulate the coupled 1D-3D flow problem into the system \eqref{eq:reform-3D}-\eqref{eq:reform-1D-bc}, wherein the singularity has been removed and all variables are smooth. 

\subsection{Elliptic equations with a single line source}
\label{sec:dec-smooth-f}
In this section, we consider the elliptic equation
\begin{align}
-\Delta p = f \delta_\Lambda
\label{eq:givenf}
\end{align}
when $\Lambda$ and $\Omega$ are as illustrated in Figure \ref{fig:illustration}, and $f=f(z) \in C^1_0(\Lambda)$ is a \textit{given}, smooth line source intensity (assumed zero at the endpoints of $\Lambda$). The solution $p$ then admits a splitting into an explicit, low-regularity term $f(z) \Psi(r) G(r)$, and an implicit, high-regularity term $v$: 
\begin{align}
p =  f(z) \Psi(r) G(r) + v(r,z).
\label{eq:u-infinite}
\end{align}
Here, $G(r)$ captures the singular part of the solution, and is given by
\begin{align}
G(r)=- \frac{1}{2\pi} \ln(r),
\end{align} 
and $\Psi(r)$ denotes some smooth cut-off function satisfying
\begin{subequations}
\begin{align}
\Psi(r) =& \, 1 & \text{ for } 0 \leq r < R_\epsilon, \\
\Psi(r) \in&  \, (0,1) \,  & \text{ for } R_\epsilon<r<R_c, \\
\Psi(r) =& \, 0 & \text{ for } r>R_c.
\end{align}
\end{subequations}
Assuming the cut-off radius $R_c$ is chosen small enough to satisfy $\Psi(r)=0$ on $\partial \Omega$, the regular component $v$ can then be defined as the solution of
\begin{subequations}
\begin{align} 
  - \Delta v &= F &  \quad \text{in } \Omega,\label{eq:w-smooth}\\ 
   v &= p_D & \quad \text{on } \partial \Omega,
\label{eq:w-smooth-bc}
\end{align}
\end{subequations}
where
\begin{equation}
F = f''(z) G(r). 
\end{equation}
To see that $p$ given by \eqref{eq:u-infinite} indeed solves \eqref{eq:givenf}, let us first note that $G = -\sfrac{1}{2\pi} \ln(r)$ was so chosen because it satisfies $-\Delta G = \delta_\Lambda$. To be more precise, $G$ is the fundamental solution of the Laplace equation in 2D, and  thus has the property
\begin{align}
-\int_{\Omega} \Delta G(r) \phi \mathrm{d} \Omega= \int_{\Lambda} \phi \mathrm{d} \Lambda \quad \forall \phi \in C^0(\Omega). \label{eq:pointsource}
\end{align}
Considering then the Laplacian of $p$ given by  \eqref{eq:u-infinite}, a straightforward calculation shows that all but one term vanish by construction, i.e.,
\begin{align}
- \Delta p &= \int_\Omega f(z) \Psi(r) \Delta G(r)  \phi \, \mathrm{d}  \Omega. 
\end{align}
By \eqref{eq:pointsource}, we then find that
\begin{align*}
- \Delta p &= \int_\Lambda f \phi \, \mathrm{d} \Lambda \quad \forall \phi \in C^0(\Omega),
\end{align*}
and it follows that the $p$ constructed in \eqref{eq:u-infinite} indeed solves \eqref{eq:givenf} in a suitably weak sense. 

Formally speaking, the splitting works by introducing first the logarithmic term $G$ for which the Laplacian returns the line source with the required intensity $f$. The higher-regularity term $v$ is then used to correct the solution so it solves the original problem. The existence of such a function $v$ follows from standard elliptic theory. As $\ln(r)\in L^2(\Omega)$, and $f''(z) \in L^2(\Lambda)$ by assumption, one can show that the entire right-hand side $F$ in \eqref{eq:w-smooth} belongs to $L^2(\Omega)$ \cite[Section 3.1]{Gjerde2018}. Consequently, there exists $v \in H^2(\Omega)$ solving \eqref{eq:w-smooth}-\eqref{eq:w-smooth-bc}. The full solution $p$, meanwhile, fails to belong to $H^1(\Omega)$. This can be shown by straightforward calculation, as one has $\ln(r) \in L^2(\Omega)$ but $\nabla \ln(r) \, \cancel{\in} \, L^2(\Omega)$. Instead, one has $p$ belonging to the weighted Sobolev space $H^1_\alpha(\Omega)$ for any $\alpha>0$. It follows that $v$ is indeed the higher-regularity term in the splitting \eqref{eq:u-infinite}. Formally, this means that $v$ is smoother and better behaved than the full solution $p$. This observation will be central to the numerical method considered in Section \ref{sec:disc}. 

\subsection{Elliptic equations with an arbitrary line source}
\label{sec:dec-smooth-f2}
In this section, we consider the elliptic problem
\begin{align}
\nabla \cdot \left( -\frac{\kappa}{\mu} \nabla p \right) = f \delta_\Lambda,
\label{eq:givenf2}
\end{align}
when the right-hand side is a line source $\delta_\Lambda$ located on a single line segment $\Lambda$ with endpoints $\mathbf{a}, \mathbf{b} \in \Omega$. The line $\Lambda$ can be described by the parametrization $\mathbf{y} = \mathbf{a} + \pmb{\tau} s \quad \text{for } s  \in (0, L),$
where $L = \Vert \mathbf{b}-\mathbf{a} \Vert$ denotes the Euclidean norm and $\pmb{\tau} = (\mathbf{b} - \mathbf{a}) / L$ is the normalized tangent vector of $\Lambda$. Letting again $f = f(s) \in C^1(\Lambda)$ be a given line source intensity, the solution $p$ then admits a splitting into an explicit, low-regularity term $E(f) G(r)$, and a high-regularity component $v$:
\begin{align}
p = E(f) \Psi G + v.
\label{eq:u-infinite2} 
\end{align}
The function $G$ is now given by
\begin{align}
G(\mathbf{x}) &=  \frac{1}{4\pi} \frac{\mu}{\kappa} \ln \left(    \frac{r_{b} + L +  \pmb{\tau} \cdot (\mathbf{a}-\mathbf{x})   } {r_{a} + \pmb{\tau} \cdot (\mathbf{a}-\mathbf{x})   }  \right),
\label{eq:G}
\end{align}
with $r_{b}(\mathbf{x}) =  \Vert \mathbf{x} - \mathbf{b} \Vert$ and $r_{a}(\mathbf{x})  =\Vert \mathbf{x} - \mathbf{a} \Vert$. This function was constructed by integrating the 3D Green's function for \eqref{eq:3D-eq} (when posed in $\mathbb{R}^2$) over the line segment $\Lambda$. It thus satisfies the property $\nabla \cdot (-\frac{\kappa}{\mu} \nabla G) = \delta_\Lambda$ \cite[Section 3.2]{Gjerde2018}. Next, $E$ denotes an extension operator $E \colon H^2(\Lambda)\rightarrow H^2(\Omega)$ extending $f$ so that it can be evaluated in the entire domain $\Omega$. Assuming again that the cut-off function $\Psi$ satisfies $\Psi=0$ on $\partial \Omega$, the regular component $v$ is then defined as the solution of 
\begin{subequations}
\begin{align} 
  - \Delta v &= F      &  \quad \text{in } \Omega,\label{eq:w-smooth2}\\ 
   v &= p_D & \quad \text{on } \partial \Omega,
\label{eq:w-smooth2-bc}
\end{align}
\end{subequations}
where 
\begin{equation}
\begin{aligned}
F &=   G \Delta \big( E(f) \Psi \big) + 2 \nabla \big( E(f) \Psi \big) \cdot \nabla G .
\end{aligned} \label{eq:FF}
\end{equation} 
To see that the constructed $p$ indeed solves the right problem, let us start by inserting it into \eqref{eq:givenf2}. construction, all terms disappear except $E(f) \Psi  \Delta G$. Integrating this term over the domain, we find that
\begin{equation}
\begin{aligned}
- \Delta p &= -\int_\Omega E(f) \Psi \Delta G \phi \, \mathrm{d}  \Omega \\ &= \int_\Lambda f \phi \, \mathrm{d} \Lambda,
\end{aligned}
\end{equation}
for all $\forall \phi \in C^0(\Omega)$, where we used the property that $E(f)=f$ on $\Lambda$. It follows that the $p$ constructed in \eqref{eq:u-infinite} indeed solves \eqref{eq:givenf} in a suitably weak sense. 

By a similar argument as the one given in \cite[Section 3.2]{Gjerde2018}, one finds that $F$ given by \eqref{eq:FF} belongs to $L^{2-\epsilon}(\Omega)$ for arbitrarily small $\epsilon>0$. It follows that there exists $v \in H^{2-\epsilon}(\Omega)$ solving  \eqref{eq:w-smooth2}-\eqref{eq:w-smooth2-bc}. Moreover, a straightforward calculation shows that $G$ again fails to belong to $H^{1}(\Omega)$. In fact, one has $G \in H^{1-\epsilon}(\Omega)$. It follows that $v$ constitutes the higher-regularity component of the solution split \eqref{eq:u-infinite2}, meaning that $v$ is smoother and better behaved than the full solution $p$.

\subsection{The coupled 1D-3D flow problem}
\label{sec:dec-coupled}
Let us now consider the coupled 1D-3D flow problem \eqref{eq:3D-eq}-\eqref{eq:1D-eq-bc}. To start with, let us again consider a single line segment $\Lambda$ with endpoints $\mathbf{a}, \mathbf{b} \in \Omega$. From the discussion in the preceding section, it is natural to assume $p$ solving \eqref{eq:3D-eq}-\eqref{eq:1D-eq-bc} admits a solution splitting of the type:
\begin{align}
p =  \Psi E(\beta f) G + v,
\label{eq:u-coupled} 
\end{align}
with $G$ being as in \eqref{eq:G}, $\Psi$ being some smooth cut-off function, $f$ being the previously introduced pressure difference $f=\hat{p} - \bar{p}$, and $v$ defined as the solution of
\begin{subequations}
\begin{align} 
 - \Delta v &= F(\hat{p}, \bar{p}; \beta) &  \quad \text{in } \Omega, \\ 
   v &=  p_D. & \quad \text{on } \partial \Omega,
\label{eq:w}
\end{align}
\end{subequations}
with
\begin{equation}
\begin{aligned}
F &=  G \Delta \big( E(\beta f) \Psi \big) + 2 \nabla \big( E(\beta f) \Psi \big) \cdot \nabla G .
\end{aligned}
\end{equation} 
Unlike in Sections \ref{sec:dec-smooth-f} and \ref{sec:dec-smooth-f2}, $f=f(\hat{p}, \bar{p})$ is now implicitly given from $\hat{p}$ and $\bar{p}$ solving the coupled 1D-3D flow problem. To reformulate \eqref{eq:3D-eq}-\eqref{eq:1D-eq-bc} in terms of $\hat{p}$ and $v$, the right-hand side therefore needs to be reformulated. To this end, let us first treat the pressure difference $\hat{p}-\bar{p}$. By the splitting \eqref{eq:u-coupled} and the definition of the averaging in \eqref{eq:avg}, calculations reveal that
\begin{equation}
\begin{aligned}
\bar{p} &=   \beta \left(\hat{p}-\bar{p}  \right) \bar{G} + \bar{v} , \\
\Rightarrow \bar{p} &= \frac{\beta \bar{G} \hat{p} + \bar{v}}{1 +  \beta \bar{G}}, \\
\Rightarrow  \hat{p} - \bar{p} &= \frac{\hat{p} - \bar{v}}{1 + \beta \bar{G}}.
\end{aligned}
\end{equation}
Here we used the simplifications $\overline{E(f)} = f \vert_\Lambda$ and $\overline{\Psi} \approx 1 \vert_\Lambda$. This is motivated by the fact that the well radius $R$ is assumed negligible. From this, we can state the reformulated coupled 1D-3D flow model: 
\begin{subequations} 
\begin{align}
- \Delta v =&   F(\hat{p}, \bar{v}; \beta^*) & \text{in} \, \Omega, \\
 v =&  p_D & \text{on} \, \partial \Omega, \\
 -\frac{\mathrm{d}^2 \hat{p}}{\mathrm{d} s^2}  =&  \hat{\beta}^* (\hat{p}-\bar{v}) & \text{in} \, \Lambda,  \\
 \hat{p} =& \hat{p}_D & \text{on} \, \partial \Lambda, 
\end{align}
\end{subequations}
where
\begin{equation}
\begin{aligned}
F &=  G \Delta \left( E(\beta^*(  \hat{p}-\bar{v})) \Psi \right)  \\ & \quad  + 2 \nabla (E(\beta^*(  \hat{p}-\bar{v})) \Psi) \cdot \nabla G ,
\end{aligned}
\end{equation}
$\beta^*$ is given by
\begin{align}
\beta^* = \frac{\beta}{1 + \beta G(R)},
\end{align} 
and $\hat{\beta}^*= \beta^* / \pi R^2$.

The extension to multiple wells follows naturally by applying the superposition principle. Considering now $\Lambda = \cup_{w=1}^{\text{wells}} \Lambda_w$, with each line segment $\Lambda_w$ having endpoints $(\mathbf{a}_w, \mathbf{b}_w) \in \Omega$, we can formulate a solution splitting 
\begin{align}
p =\sum_{w=1}^{\text{wells}} E \big( \beta^* (\hat{p}-\bar{v}) \big) \Psi_w G_w + v ,
\label{eq:u-existence}
\end{align} 
where $E \colon H^2(\Lambda) \rightarrow H^2(\Omega)$ is the same extension operator as before,$G_w$ is given by \eqref{eq:G} with $\mathbf{a}=\mathbf{a}_w$ and $\mathbf{b}=\mathbf{b}_w$, $\Psi_w$ is some smooth cut-off function with respect to line segment $\Lambda_w$, and $v$ solves
\begin{subequations}
\begin{align}
- \Delta v =&  F(\hat{p}, \bar{v}; \beta^*) & \text{in} \, \Omega,\label{eq:reform-3D} \\
 v =& p_D & \text{on} \, \partial \Omega,\label{eq:reform-3D-bc}\\
- \frac{\mathrm{d}^2 \hat{p}}{\mathrm{d} s^2}  =& -  \hat{\beta}^{*} ( \hat{p} - \bar{v}) & \text{in} \, \Lambda,\label{eq:reform-1D} \\
 \hat{p} =& \hat{p}_D & \text{on} \, \partial \Lambda,\label{eq:reform-1D-bc}
\end{align}
\end{subequations}
with right-hand side 
\begin{equation}
\begin{aligned}
F &=\sum_{w=1}^\text{wells}   G_w \Delta \big( E(\beta^*(  \hat{p}-\bar{v})) \Psi_w \big)  \\ & \quad  + 2 \nabla \big( E(\beta^*(  \hat{p}-\bar{v})) \Psi_w \big) \cdot \nabla G_w ,
\end{aligned}
\end{equation}
and 
\begin{align}
\beta^* = \frac{\beta}{1 + \sum_{w=1}^{\text{wells}}  \beta \overline{G_w \Psi_w}}, \quad \hat{\beta}^* = \frac{\beta^*}{\pi R^2}.
\end{align}
The system \eqref{eq:reform-3D}-\eqref{eq:reform-1D-bc} constitutes a reformulation of the coupled 1D-3D flow model in terms of the smooth variables $v$ and $\hat{p}$. For an example of what the splitting might look like, the reader is invited to examine Figure \ref{fig:well1}. As the singularities have here been removed from the system, it enjoys significantly improved regularity compared to the standard formulation \eqref{eq:3D-eq}-\eqref{eq:1D-eq-bc}.

\section{Weak formulation}
\label{sec:variational-form}
In this section, we state a weak formulation of the reformulated coupled 1D-3D flow problem \eqref{eq:reform-3D}-\eqref{eq:reform-1D-bc}. As the variables in this formulation are all smooth functions, this can be done using standard Sobolev spaces. For the sake of completeness, we give also a weak formulation of the standard coupled 1D-3D flow problem \eqref{eq:3D-eq}-\eqref{eq:1D-eq-bc}. The reservoir pressure $p$ therein contains a singularity; for this reason, its weak formulation requires the use of weighted Sobolev spaces.

Consider first the reformulated coupled 1D-3D flow problem.
Let $\mathbf{V}$ denote the product space $\mathbf{V} = V \times \hat{V}$, where
\begin{align}
V &= \lbrace u \in H^1(\Omega) \colon u \vert_{\partial \Omega} =  p_D \rbrace, \\
\hat{V} &= \lbrace \hat{u} \in H^1(\Lambda) \colon u \vert_{\partial \Lambda} =  \hat{p}_D \rbrace,
\end{align}
normed by
\begin{align}
\norm{(v, \hat{p})}{\mathbf{V}}^2 = \norm{v}{H^1(\Omega)}^2 + \norm{\hat{p}}{H^1(\Lambda)}^2.
\end{align}

Multiplying \eqref{eq:reform-3D} and \eqref{eq:reform-1D} with test functions $\phi \in H^1_0(\Omega)$ and $\hat{\phi} \in H^1_0(\Lambda)$, respectively, integrating over their respective domains, and performing an integration by parts, we arrive at the following variational formulation:

Find $(v, \hat{p}) \in \mathbf{V}$ such that
\begin{align}
a\left( \left( v, \hat{p}  \right), ( \phi, \hat{\phi}  ) \right) = 0
\label{eq:varform-a}
\end{align}
for all $(\phi, \hat{\phi}) \in \mathbf{V}_0$, where 
\begin{equation}
\begin{aligned}
a  \left( \left( v, \hat{p}  \right),  \phi, \hat{\phi}  ) \right) &= \left( \nabla v, \nabla \phi \right)_\Omega + \left( \frac{\mathrm{d}}{\mathrm{d} s} \hat{p}, \frac{\mathrm{d}}{\mathrm{d} s} \hat{\phi}\right)_\Lambda \\
  &+  \left( F_1( \beta^*(\hat{p} - \bar{v})) , \nabla \phi\right)_\Omega  \\
  &-  \left(  F_2(\beta^*(\hat{p} - \bar{v})) , \phi \right)_\Omega  \\
  &+ (\hat{\beta}^{*} (\hat{p}-\bar{v})  ), \hat{\phi})_\Lambda ,
\end{aligned}
\end{equation}
and 
\begin{subequations}
\begin{align}
F_1(\hat{\phi}) &= \sum_{w=1}^\text{wells} \nabla \Big( \Psi_w E \big( \hat{\phi} \big) \Big) G_w, \\
F_2(\hat{\phi}) &= \sum_{w=1}^\text{wells} \nabla \Big( \Psi_w E \big( \hat{\phi}) \big)\Big)  \cdot \nabla G_w.
\end{align} 
\end{subequations} 
The full reservoir pressure can then be constructed from $v$ and $\hat{p}$ by the relation
\begin{align}
p =  \sum_{w=1}^{\text{wells}} E \big( \beta^* \left( \hat{p}-\bar{v}\right) \big) G_w + v.
\label{eq:u-full-discrete}
\end{align}

Next, let us consider the standard coupled 1D-3D flow model, and give its variational formulation as it was proposed in \cite{dangelo2008}. Let $\mathbf{V}_\alpha$ denote the weighted product space $\mathbf{V}_\alpha = V_\alpha \times \hat{V}$,
where 
\begin{align}
V_\alpha &= \lbrace u \in H^1_\alpha(\Omega) \colon u \vert_{\partial \Omega} = p_D \rbrace, \\
\hat{V} &= \lbrace \hat{u} \in H^1(\Lambda) \colon \hat{u} \vert_{\partial \Lambda} = \hat{p}_D \rbrace, 
\end{align}
normed by 
\begin{align}
\norm{(p, \hat{p})}{\mathbf{V}_\alpha}^2 = \norm{p}{H^1_\alpha(\Omega)}^2 + \norm{\hat{p}}{H^1(\Lambda)}^2.
\end{align}
Multiplying \eqref{eq:3D-eq} and \eqref{eq:1D-eq} with test functions $v \in H^1_{-\alpha, 0}(\Omega)$ and $\hat{v} \in H^1_0(\Lambda)$, respectively, integrating over their domain of support, and performing an integration by parts, we arrive at the variational formulation:

Find $(p, \hat{p}) \in \mathbf{V}_\alpha$ such that
\begin{align}
a\left( ( p, \hat{p} ), ( \phi, \hat{\phi}  ) \right) = 0
\label{eq:varform-a-standard}
\end{align}
for all $(\phi, \hat{\phi}) \in \mathbf{V}_{-\alpha, 0}$, where 
\begin{equation}
\begin{aligned}
a\left( \left( p, \hat{p}  \right), ( \phi, \hat{\phi}  ) \right) &= \left( \nabla p, \nabla \phi\right)_\Omega + \left( \frac{\mathrm{d}}{\mathrm{d} s} \hat{p}, \frac{\mathrm{d} }{\mathrm{d} s} \hat{\phi}\right)_\Lambda \\
  &- \left( \beta \left( \hat{p}-\bar{p}\right), \phi\right)_\Lambda \\
  &+  (\hat{\beta} (\hat{p}-\bar{p})  ), \hat{\phi})_\Lambda,\label{eq:a-standard}
\end{aligned}
\end{equation}
and the test space $\mathbf{V}_{-\alpha, 0}$ is the space of functions $(\phi, \hat{\phi})\in \mathbf{V}_{ -\alpha, 0}$ with zero trace on the boundary. Notice here that the test and trial spaces are chosen with opposite weight functions; this is what ensures the continuity and coercivity of the bilinear form \eqref{eq:a-standard}. For a proof of the well-posedness of this formulation, the reader is referred to \cite{dangelo2008, dangelo2012}.

\section{Numerical Discretization}
\label{sec:disc}
In this section, we show the block matrix resulting from a finite element discretization of weak formulation of the reformulated coupled 1D-3D problem. As the pressure difference $f(\hat{p}, \bar{v})=\hat{p}-\bar{v}$ now uses the regular part of the pressure, $v \in H^2(\Omega)$, we introduce here also the simplification $\bar{v}_h=v_h \vert_\Lambda$; i.e., we take the trace of $v_h$ on $\Lambda$ rather than the average over the cylinder. This is motivated by the fact that $R$ is assumed negligible compared to the mesh size $h$, and $v$ is regular, meaning $\bar{v} \approx  v \vert_\Lambda$. The result is a \enquote{true} coupled 1D-3D flow model, in that it considers only 1D and 3D variables, with no averaging performed over a 2D cylinder. The same approximation is not possible for the standard coupled 1D-3D flow model as the reservoir pressure is there undefined on $\Lambda$.

We will now give the discretized form of the variational formulation \eqref{eq:varform-a}. For simplicity, let us assume $\Omega$ is a polyhedron that readily admits a partitioning $\mathcal{T}_{T, h}$ into simplicial elements $T$:
\begin{align*}
\bar{\Omega} = \bigcup_{T \in \mathcal{T}_{T, h}} T.
\end{align*}
The simplicial partitioning $\mathcal{T}_{T, h}$ forms a mesh, assumed conforming, which can then be characterized by the mesh size $h= \max_{T \in \mathcal{T}_{T,h}} h_T$. Next, we associate  this mesh with the usual (3D) Lagrange space of order 1, $V^{h}_{u}$, given by 
\begin{align*}
V^{h}_{u} = \lbrace v_h \in C^0_{u}(\Omega), \, v_h \vert_T \in \mathbb{P}_1  \text{ where } T \in \mathcal{T}_{T, h} \rbrace.
\end{align*}
Here,  $\mathbb{P}^1$ denotes the space of polynomials of degree $1$, and $C^0_{u}(\Omega)$ the space of continuous elements that equal the interpolation of $u$ on the boundary, i.e., 
\begin{align}
C^0_{u}(\Omega) = \lbrace p \in C^0(\Omega) \colon p \vert_{\partial \Omega} = \mathcal{I}_h u \rbrace.\label{eq:c0}
\end{align} 
Next, we assume $\Lambda$ admits a partitioning $\mathcal{T}_{I, h}$ into line segments $I$:
\begin{align*}
\bar{\Lambda} = \bigcup_{I \in \mathcal{T}_{I, h}} I,
\end{align*}
assumed again to satisfy all the requirements of a conforming mesh, and associated with the mesh size $\hat{h} = \max_{I \in \mathcal{T}_{I,h}} h_I$. For the discretization of $\hat{V}$, we use the (1D) Lagrange space of order 1, 
\begin{align*}
\hat{V}^{h}_{\hat{p}} = \lbrace v_h \in C^0_{\hat{u}}(\Lambda), \, \hat{v} \vert_I \in \hat{\mathbb{P}}_1  \text{ where } I \in \mathcal{T}_{I, h} \rbrace,
\end{align*}
with $C^0_{\hat{u}}(\Lambda)$ interpreted as in \eqref{eq:c0}. 

Considering first the reformulated system \eqref{eq:varform-a}, let  
\begin{align}
v = \sum_{k=1}^N v_k \phi_k, \qquad \hat{p} = \sum_{l=1}^{\hat{N}} \hat{p}_l \hat{\phi}_l,
\end{align}
where $\lbrace \phi_1, \phi_2, ..., \phi_N \rbrace$ and $\lbrace \hat{\phi}_1, \hat{\phi}_2, ..., \hat{\phi}_{\hat{N}} \rbrace$ are linear hat functions spanning $V^h$ and $\hat{V}^h$, respectively. Note next that $v_h$ is a linear function used to approximate the high regularity term $v \in H^2(\Omega)$. For $R \ll h$, its average $\bar{v}_h$ can be well approximated by simply taking the trace $v \vert_\Lambda$. The pressure difference $\hat{p}-\bar{v}$ is then given by
\begin{equation}
\begin{aligned}
\hat{p}-\bar{v} &= \sum_{l=1}^{\hat{N}} \hat{p}_l \hat{\phi}_l - \sum_{k=1}^N v_k \bar{\phi}_k \\
&= \sum_{l=1}^{\hat{N}} \hat{\phi}_l \left( \hat{p}_l   - \sum_{k=1}^N  T_{k,l}  v_k \right)  .
\end{aligned}
\end{equation}
Here, $T \colon V^h \rightarrow \hat{V}^h$ is the discrete trace matrix, characterized by the property $\phi_k \vert_\Lambda = \sum_{l=1}^{\hat{N}} T_{k,l} \hat{\phi}_l$. 

Testing \eqref{eq:varform-a} with $v = \phi_i$ for $i=1,...,N$ and $\hat{v} = \hat{\phi}_j$ for $j=1,..., \hat{N}$, we arrive at the following discrete system: 
\begin{align}
\begin{bmatrix}
  A - CT^T  & C\\
 -\hat{M} T^T   & \hat{A} + \hat{M}
\end{bmatrix}
\begin{bmatrix}
   v \\
   \hat{p}
\end{bmatrix}
= 0.\label{eq:block2}
\end{align} 
where $A$ and $\hat{A}$ are the standard stiffness matrices
\begin{align}
A_{i,k} &= ( \nabla \phi_k , \nabla \phi_i), \\
\hat{A}_{j,l} &= ( \frac{\mathrm{d}}{\mathrm{d} s} \hat{\phi}_l, \frac{\mathrm{d}}{\mathrm{d} s} \hat{\phi}_j).
\end{align}
$\hat{M}$ denotes the standard 1D mass matrix,
\begin{align}
\hat{M}_{j,l} = ( \hat{\beta}^* \hat{\phi}_j, \hat{\phi}_l)_\Lambda
\end{align}
and $C$ denotes the coupling block, 
\begin{align}
C_{i, l}&= \sum_{w=1}^\text{wells}  \left( F_1(\beta^* \hat{\phi}_l) , \nabla \phi_i \right)_\Omega - \left( F_2 ( \beta^* \hat{\phi}_l ) , \phi_i \right)_\Omega.
\end{align}

We will refer to this system as the Singularity Removal Based FE method. After solving \eqref{eq:block2}, a discretization of the full reservoir pressure $p_h$ can be reconstructed using 
\begin{align}
p_h  = \sum_{w=1}^{\text{wells}} \beta^* ( \hat{p}_h- \bar{v}_h \vert_\Lambda) \mathcal{I}_h^k G  + v_h  \label{eq:p-reconstruction},
\end{align}
where $\mathcal{I}_h^k$ denotes the interpolation onto the Lagrange space of order $k$. As the interpolation of $G(r)$ is fairly cheap, the approximation property of $p_h$ can here be improved by choosing the interpolation degree $k$ high.

A more straightforward method can be found by discretizing \eqref{eq:varform-a-standard} directly; this is the finite element formulation analysed in e.g. \cite{dangelo2012}. As we will compare the performance of this method against the Singularity Removal Based FE method, we give here its discretization for the sake of completeness. Setting
\begin{align}
p = \sum_{k=1}^N p_k \phi_k \qquad \hat{p} = \sum_{l=1}^{\hat{N}} \hat{p}_l \hat{\phi}_l,
\end{align}
The pressure difference $\hat{p}-\bar{p}$ is then given by
\begin{equation}
\begin{aligned}
f &= \hat{p} - \bar{p} \\
&= \sum_{l=1}^{\hat{N}} \hat{p}_l \hat{\phi}_l - \sum_{k=1}^N p_k \bar{\phi}_k \\
&= \sum_{l=1}^{\hat{N}} \hat{p}_l \hat{\phi}_l  - \sum_{k=1}^N \sum_{m=1}^{\hat{M}} \Pi_{m,k} p_k   \hat{\psi}_m,
\end{aligned}
\end{equation}
where $\Pi$ is the discrete averaging matrix $\Pi : V^h \rightarrow X^h$ and $\{ \hat{\psi}_1, \hat{\psi}_2, ..., \hat{\psi}_{\hat{M}} \}$ are the basis functions spanning $\hat{X}^h$. 

Testing now \eqref{eq:varform-a-standard} with $v = \phi_i$ for $i=1,...,N$ and $\hat{v} = \hat{\phi}_j$ for $j=1,..., \hat{N}$, we arrive at the following block system for the discretization of \eqref{eq:varform-a-standard}:
\begin{align}
\begin{bmatrix}
  A + \beta T^T N \Pi  \quad & -\beta T^T \hat{M}\\
 - \beta N \Pi  & \hat{A} + \beta \hat{M}
\end{bmatrix}
\begin{bmatrix}
   p \\
   \hat{p}
\end{bmatrix}
= 0.\label{eq:block1}
\end{align}
Here, $N$ denotes the mass matrix given by 
\begin{align}
N_{m, l} = (\hat{\psi}_m, \hat{\phi}_l),
\end{align}
for $\hat{\psi}_m$ belonging to the discontinuous Galerkin space of order 0:
\begin{align*}
\hat{X}^{h} = \lbrace v_h \in L^2(\Gamma), \, v_h \vert_I \in \mathbb{P}_0  \text{ where } I \in \mathcal{T}_{I, h} \rbrace.
\end{align*}
We will refer to this system as the standard FE method.

\section{Relation to the Peaceman well model}
\label{sec:peaceman}
In this section, we show that the reformulated coupled 1D-3D flow model \eqref{eq:reform-3D}-\eqref{eq:reform-1D-bc} under certain conditions reduces to the Peaceman well correction. We start by giving a brief summary of the methodology Peaceman introduced in his seminal work \cite{peaceman1978}. We then return to our reformulated model, and show that with $G(r)$ chosen so that its support is the equivalent radius of the Peaceman well correction, the reformulation results in a well index that equals the one derived by Peaceman.

In reservoir simulations, the mass flux between well and aquifer,  $q$, is usually modelled in a manner analogous to that in \eqref{eq:3D-eq}:
\begin{align}
q = J (p_w-p_K),
\label{eq:wellflow}
\end{align}
where $p_{w}$ is the flowing pressure in the well, $J$ its well index, and $p_{K}$ the reservoir pressure averaged over the grid cell $K$. In Section \ref{sec:decomposition}, we showed how the line source that models the well introduces a logarithmic type singularity in the reservoir pressure. For wells with radius much smaller than the grid size $h$, i.e., $R \ll h$, $p_{K}$ is therefore likely to constitute a poor representation of the reservoir pressure in the near vicinity of the well. 

The Peaceman well model accounts for this by altering the well index $J$ in \eqref{eq:wellflow} so that $q$ better corresponds to the numerical approximation of the pressure difference between well and aquifer. Assuming radial flow, Darcy's law in a heterogeneous reservoir is given, per unit well length, by the relation
\begin{align}
\frac{q}{2\pi r} = - \frac{\kappa}{\mu} \frac{\mathrm{d}p}{\mathrm{d}r}.
\end{align}
Integrating this equation to a radius $r_e$,
\begin{align}
 \frac{2\pi \kappa}{q \mu} \int_{p_{w}}^{p_e} \mathrm{d}p =- \int_{R}^{r_e} \mathrm{d}r,
\end{align}
we find that
\begin{align}
q = \frac{2\pi \kappa}{\mu} \frac{ p_{w}-p_e}{\ln(r_e/R)}
\end{align}
when $p_e=p(r_e)$. We also need to take into account the pressure drop $\Delta p_{\text{skin}}$ across the skin of the well. To do so, let $S$ be the skin-factor, defined by the relation
\begin{align}
S = \frac{2\pi \kappa}{q \mu} \Delta p_{\text{skin}}.
\end{align}
Letting now $r_e$ be the radius at which the reservoir pressure equals the averaged grid cell pressure $p_K$, Peaceman used the following relation between $q$ and the pressure difference $p_{w}-p_K$ \cite{peaceman1978}:
\begin{align}
q = \frac{2\pi \kappa}{\mu} \frac{ p_{w}-p_K }{\ln(r_e/R) + S}.
\label{eq:peace}
\end{align}

To utilize this correction, one must first identify the equivalent radius $r_e$ entering in \eqref{eq:peace}. This radius generally depends on the discretization method, the location of the well within the grid, and the permeability of the rock around the well. Assuming for example square grid blocks and a well at the center of an interior grid block, Peaceman derived an equivalent radius $r_e = 0.2h$ for the two-point flux approximation \cite{peaceman1978}.

\begin{figure*}
\centering
\includegraphics[width=0.9\textwidth]{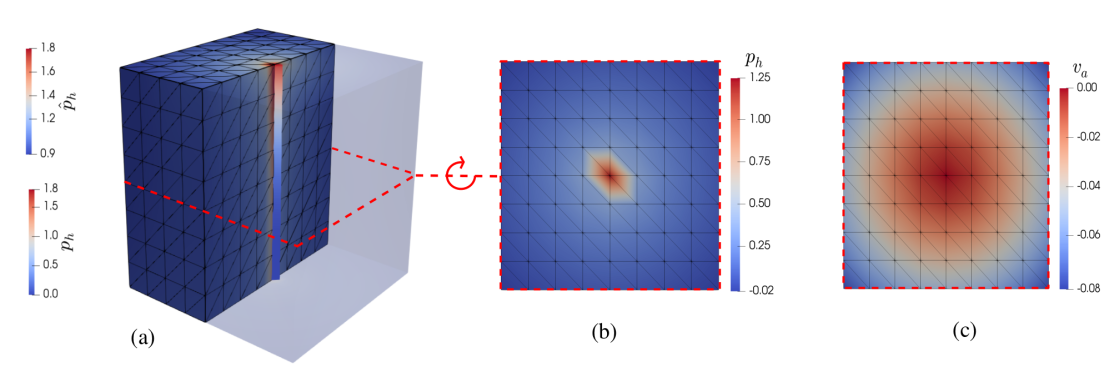}
\caption{(a) FE approximations of $\hat{p}_h$ and the reconstructed reservoir pressure $p_h$ for $h = \sfrac{1}{8}$. (b) Full reservoir pressure $p_h$ and (c) background pressure $v_h$ on the slice $\{ (x,y,z) \in \Omega : z = 0.5\}$}
\label{fig:well1}
\end{figure*}

The reformulation of the pressure difference $f$ in terms of $\hat{p}$ and $v$ bears a strong resemblance to the Peaceman well correction in \eqref{eq:peace}. In a practical sense, the reformulation into \eqref{eq:reform-3D}-\eqref{eq:reform-1D-bc} can be interpreted as a non-local well correction, which has a support in a region around the well which may significantly exceed the grid resolution. To see more clearly the similarity with the Peaceman well correction, let us now consider a single well. We have then
\begin{align}
q = \frac{\beta}{1 + \frac{\mu}{\kappa}\beta \overline{G}} (\hat{p}-\bar{v}).
\label{eq:mypeace}
\end{align}
Next, we let now $\hat{p}$ be the flowing well pressure $p_{w}$. The term $G(r)$ contains the logarithmic component of the solution; in a manner analogous to the Peaceman well correction, we make it local to the cylinder of radius $r_e$ by setting
\begin{align}
G_{r_e}(r) = \begin{cases}
-\frac{1}{2\pi} \frac{\mu}{\kappa} \ln(r/r_e) &\text{ for } r \leq r_e,\\
0 &\text{ otherwise}.
\end{cases}  
\end{align}
Note that this $G$ is not smooth enough to work for the solution split \eqref{eq:u-existence}, we use it here only for the sake of comparison. By the definition of the averaging \eqref{eq:avg}, we have $\overline{G} = - \sfrac{\mu}{2\pi \kappa} \ln(\sfrac{R}{r_e})$ Inserting it in \eqref{eq:mypeace} yields the relation
\begin{align}
q
&= \frac{\beta}{1 - \beta \frac{\mu}{2\pi \kappa} \ln(\frac{R}{r_e}) } (p_{w}-\bar{v}) \\
&= \frac{2 \pi \kappa}{\mu} \frac{p_{w}-\bar{v}}{\frac{2\pi \kappa }{\mu\beta} + \ln(\frac{r_e}{R})}.
\end{align}
Here, $\sfrac{2\pi \kappa}{\mu \beta}$ can be substituted by the skin factor of the well by recalling $q = \beta \Delta p_{\text{skin}}$. This results in an expression that equals the Peaceman well correction given in \eqref{eq:peace}, i.e., 
\begin{align}
q &= \frac{2 \pi \kappa}{\mu} \frac{p_{w}-\bar{v}}{\ln(\frac{r_e}{R}) + S}.
\end{align}
The regular component $v$ can be interpreted as a sort of background pressure, or more precisely, the component of the reservoir pressure that can be approximated using linear functions. We see then that the singularity removal constitutes an alteration of $\beta$ (which can be interpreted as a well index) so that the mass flux function $q$ better corresponds to the \textit{numerically} computed pressure difference between well and reservoir, i.e., $\hat{p}-\bar{v}$. For this reason, we expect that the singularity removal, in a manner similar to the Peaceman well correction, will improve the stability of the FE approximation with respect to the ratio $R/h$.

\section{Numerical Results}
\label{sec:numerics}

\label{sec:numerics1}
In this section, we perform numerical experiments to test the approximation properties of the Singularity Subtraction Based FE method given by \eqref{eq:block2}. For the implementation, we utilized the finite element framework FEniCS \cite{LoggMardalEtAl2012a}. For the first test case, we consider a single well with smooth lateral well permeability $\beta$, and compare the results against those obtained using the standard FE given by \eqref{eq:block1}. Our implementation of this method uses an earlier implementation from Kuchta \cite{miro-git}, the same as was utilized for the results of Holter et al. in \cite{Holter2018}. The Singularity Removal Based FE method was implemented by an extension of this code, using also the mixed-dimensional functionality of FEniCS developed and implemented by Daversin-Catty \cite{cecile-git}. For the second test case, we consider a discontinuous lateral permeability $\beta$, and an extension operator that uses radial basis function interpolation. We show here that the reconstructed reservoir pressure $p_h$ converges optimally when the Singularity Removal Based FE method is applied.

\subsection{Convergence test for well with smooth lateral permeability}

In this section, we take $\Omega=(0,1)^3$ and $\Lambda = \{ (x,y,z) \in \Omega : x=y= \sfrac{1}{2} \}$. We want to test the capability of each method in approximating the test problem
\begin{subequations}
\begin{align}
p_a &= -\frac{1}{2\pi}  (z^3+1) \ln(r) + v_a ,\label{eq:testproblem-p}\\
v_a &= -\frac{3}{4 \pi} \left( z r^2 \left( \ln(r)-1 \right) \right),\label{eq:testproblem-w}\\
\hat{p}_a &= \frac{1-\ln(R)}{2\pi} \left( z^3 + 1 -\frac{3}{2}R^2z \right),\label{eq:testproblem-hat}
\end{align}
\end{subequations}
with the following parameters:
\begin{align}
\kappa = \hat{\kappa}= \mu = 1, \quad \beta = {2\pi}, \quad \hat{\beta}=\frac{6z(1-\ln(R))}{z^3+1}.
\end{align}
The solution, along with the splitting terms, are shown in Figure \ref{fig:well1}.

In order to test the stability of the approximation when the well radius is small compared to mesh size $h$, we test using four different values for the well radius:
\begin{align}
R\in \{ 1.0\text{e-}1, \,  1.0\text{e-}2, \, 1.0\text{e-}3,  \, 1.0\text{e-}4 \}.
\end{align}
Furthermore, we set $\Psi = 1$ and choose as the extension operator
\begin{align}
E(f) = f(z) \text{ for all } (x,y,z) \in \Omega.
\end{align}
In this case, the reformulated FE method will approximate the analytic solution for $v_a$ given in \eqref{eq:testproblem-w}, meaning we can compute its error directly using $\norm{v_a-v_h}{}$.

\begin{figure*}[t]
\centering
\begin{subfigure}{0.34\textwidth}
\includegraphics[width=0.9\textwidth]{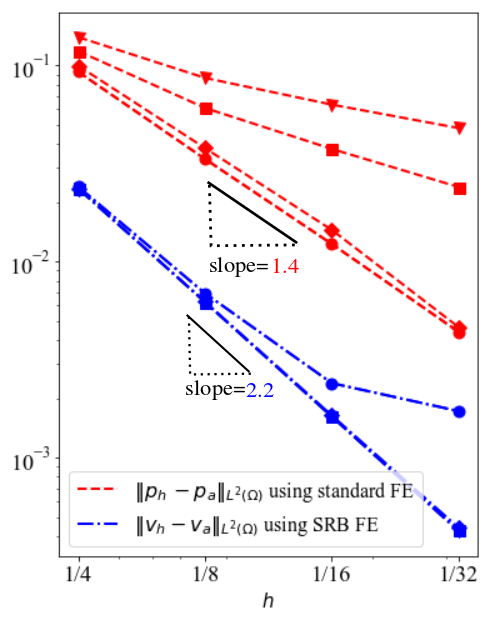}
\caption{Reservoir pressure approximations, $p-p_a$ and $v-v_a$.}\label{fig:res-error}
\end{subfigure} \hspace{0.5em}
\begin{subfigure}{0.36\textwidth}\includegraphics[width=0.9\textwidth]{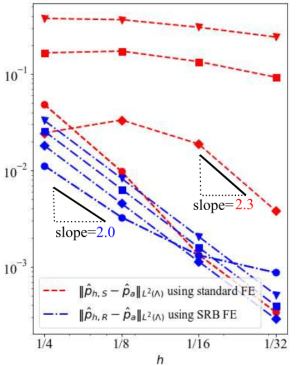}
\caption{Well pressure approximation errors, $\hat{p}_h-\hat{p}_a$.}\label{fig:well-error}
\end{subfigure} 
\begin{subfigure}{0.07\textwidth}\includegraphics[width=0.99\textwidth]{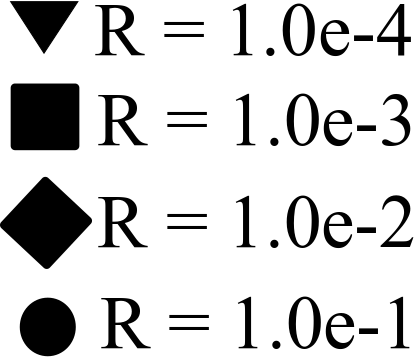}
\end{subfigure}

\caption{Log-log plot of the approximation errors obtained using the standard FE method (red) and the Singularity Removal Based FE method (blue) as the mesh size $h$ decreases. The approximations were tested for different well radius values $R$ and is indicated with a marker, where the radius corresponding to each marker is shown to the right.}\label{fig:errors}
\end{figure*}

Figure \ref{fig:errors} shows the approximation errors, measured in the $L^2$-norm, when the problem was solved using a sequence of increasingly fine meshes. The blue lines in Figure \ref{fig:res-error} show the approximation error of $v_h$, measured in the $L^2$-norm, i.e., $\norm{v_h-v_a}{L^2(\Omega)}$ with $v_a$ being the analytic solution in \eqref{eq:testproblem-w}. For $R<h$, the errors are seen to be invariant with respect to $R$, and the approximation of $v_h$ exhibits moderate superconvergence. To expand upon this, we expect for this approximation optimal convergence rates of order $h^l$ with $l=2.0$; we see here a slight super-convergence as $l=2.2$. For $h>1/8$ and $R=0.1$, our assumption of $R<h$ is no longer valid, and we see a degradation of the convergence rates. To be more precise, we made in the construction of the block matrix \eqref{eq:block2} the simplification $\bar{v}=v\vert_\Lambda$, and this is not justified for $R \sim h$. Optimal convergence rates could be restored by taking the average of $v_h$ rather than its trace.

The red lines in Figure \ref{fig:res-error} give the approximation errors for the full reservoir pressure using the standard FE method described by \eqref{eq:block1}. We give here the approximation error of $p_h$ in the $L^2$-norm, i.e., $\norm{p_h-p_a}{L^2(\Omega)}$. For the standard FE method, the convergence properties strongly depend on the well radius $R$, with decreasing $R$ leading to a reduction in the convergence rate. The best convergence rates are seen when $R \sim h$, but even here, the convergence is sub-optimal compared to the Singularity Removal Based FE method. This can be explained by noting that the standard FE method explicitly resolves the line source in the problem; it was shown in \cite{dangelo2012} that this leads to a reduction in the convergence rate of $p_h$. We refer here to our comments in \cite[p. 14-15]{Gjerde2018} for a more in-depth explanation of this, and remark only that the line source problem is expected to converge with order $h^{1-\epsilon}$ for $\epsilon>0$ arbitrarily small. Thus, the convergence order $h^l$ with $l=1.4$ surpasses the theoretical expectation when $R \sim h$. 

The blue and red lines in Figure \ref{fig:well-error} give the approximation error of $\hat{p}_h$ using the Singularity Removal Based and standard FE method, respectively. The approximation error is also here measured in the $L^2$-norm, i.e., using $\norm{\hat{p}_h-\hat{p}_a}{L^2(\Lambda)}$. We see here that the singularity removal significantly improves the convergence properties of the problem for $R<h$. The convergence rates degrade when $R>h$. This is again due to the simplification $\bar{v}=v\vert_\Lambda$ used in the construction of the block matrix \eqref{eq:block2}, and could be resolved by removing this simplification.

From \ref{fig:well-error}, is clear that the standard FE method has trouble approximating the solution when $R < h$. Moreover, the approximation error of $\hat{p}_h$ is seemingly more sensitive than $p_h$ with respect to the ratio $R/h$. This can be understood by returning to the reservoir pressure splitting $p =  \beta (\hat{p}-\bar{p}) G + v$, where $G = - \sfrac{1}{2\pi}\ln(r)$, and noting that the error in $p_h$ is due to three separate issues, namely, the error in the approximation of the pressure difference, i.e., $\norm{\hat{p}_h-\bar{p}_h-(\hat{p}_a-\bar{p}_a)}{L^2(\Lambda)}$, the error in approximating the logarithm, i.e., $\norm{\ln(r)_h- \ln(r)}{L^2(\Omega)}$, and the error in approximating $v$ (which is comparatively small). The standard FE method has trouble resolving the logarithmic nature of the reservoir pressure around the well, leading to a large approximation error in $\bar{p}$. This further pollutes the approximations of both $\hat{p}$ and $p$. The effect is not as noticeable for $p$ as its approximation error is dominated by the approximation error for $\ln(r)$. The well pressure $\hat{p}$, however, is in principle a smooth function, for which the FE approximation should  be comparatively small. Its approximation error is therefore dominated by the term $\norm{\bar{p}_a- \bar{p}_h}{L^2(\Lambda)}$. 

In summary, we see here that the standard FE method has difficulty resolving the pressure difference $\hat{p}- \bar{p}$ when $R< h$, due to the fact that $\bar{p}$ is then poorly approximated. This further pollutes the approximations of both the well and reservoir pressure. Explicitly subtracting the singularity in $p$, which results in the Singularity Removal Based FE described by \eqref{eq:block2}, restores optimal convergence rates for the reservoir pressure $p$, and improves the robustness of the method with respect to a small well radius $R$.

\subsection{Convergence test for well with discontinuous lateral permeability}

\begin{figure}[t]
\centering
\includegraphics[width=0.4\textwidth]{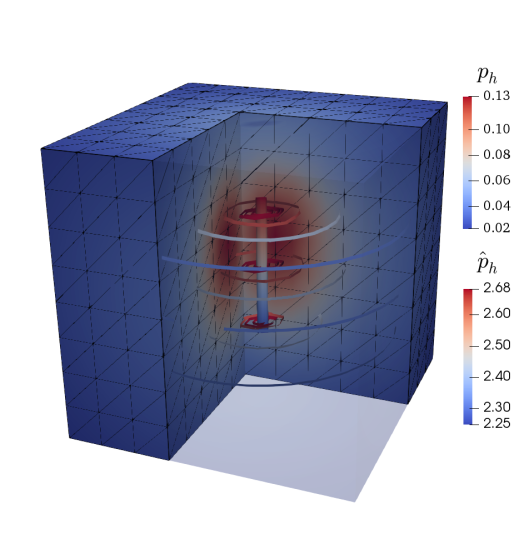}
\caption{SRB FE approximations of the reconstructed reservoir pressure $p_h$ and well pressure $\hat{p}_h$. Isolines are plotted for $p_h$.} \label{fig:exp2-pressure}
\end{figure}
Let $\Omega = (0,1)^3$ and 
\begin{align}
\Lambda = \{ (x,y,z) \in \Omega: x=y=\frac{1}{2}, \, z \in (\frac{1}{4}, \frac{3}{4}) \}.
\end{align}
In this section, we will test the ability of the Singularity Remov1l Based FE method in approximating the analytic test problem 
\begin{subequations}
\begin{align}
p_a &=  z G + v_a ,\label{eq:testproblem2-p}\\
v_a &= \frac{1}{4 \pi} (r_b-r_a),\label{eq:testproblem2-w}\\
\hat{p}_a &= \sin(z) + 2,\label{eq:testproblem2-hat}
\end{align}
\end{subequations}
when $G$ is given as in Section \ref{sec:dec-smooth-f2}:
\begin{align}
G = \frac{1}{4\pi} \ln \Big(   \frac{r_b - (z-b) }{r_a - (z-a)} \Big).
\end{align}
The problem parameters are then as follows:
\begin{align}
\kappa = \hat{\kappa}= \mu = 1, \quad \beta = \frac{z}{\hat{p}_a-\bar{p}_a}, \quad \hat{\beta}=- \frac{\beta \sin(z)}{z}.
\end{align}
Physically, this can be interpreted as modelling a well that passes through the domain but is only in contact with the reservoir when $\sfrac{1}{4} < z < \sfrac{3}{4}$. This translates to a jump in the lateral permeability, with discontinuities at the points $(\sfrac{1}{2}, \sfrac{1}{2}, \sfrac{1}{4})$ and $(\sfrac{1}{2}, \sfrac{1}{2}, \sfrac{3}{4})$. 

\begin{table}[]
\centering
\begin{tabular}{lrrrr} \normalsize
$h$ & $\norm{p_e}{L^2(\Omega)}$ & $\norm{p_e}{H^1(\Omega)}$ & $\norm{\hat{p}_e}{L^2(\Lambda)}$ & $\norm{\hat{p}_e}{H^1(\Lambda)}$ \\ \hline
$\sfrac{1}{4}$ & 1.94e-02  & 1.88e-01 & 2.30e-3 & 2.51e-2  \\
$\sfrac{1}{8}$ & 5.44-03  &  4.99e-02 & 6.27e-4 & 1.26e-2  \\
$\sfrac{1}{16}$ & 1.25e-03 & 9.26e-02 & 1.55e-4 & 6.27e-2 \\
$\sfrac{1}{32}$ & 2.77e-04 & 4.50e-02 & 7.80e-5 & 3.32e-2 \\ \hline \hline
\normalsize $l$ & 2.0 & 1.0 & 2.0 & 1.0 
\end{tabular}
\caption{The reservoir pressure approximation error $p_e = \mathcal{I}_h(p_a) - p_h$ and well pressure approximation error $\hat{p}_e = \hat{p}_a - \hat{p}_h$ when $p_h$ was reconstructed using \eqref{eq:reconstruct} with $k=1$. Both errors were found to converge with optimal order, i.e., with $l=2$ in the $L^2$-norm and $l=1$ in the $H^1$-norm.} \label{tab:exp2}
\end{table}

As the cut-off function, we use the Gaussian function 
\begin{align}
\Psi =  \exp( - \frac{\text{dist}(\mathbf{x}, \Lambda)^2 }{2c^2} )
\end{align}
with $c=0.04$. For the extension operator $E$, we choose spline interpolation with radial basis functions as given in \cite{splineinterpolation}. Given a discretized solution pair $(v_h, \hat{p}_h)$ to \eqref{eq:block2}, we can then reconstruct the discretized full reservoir pressure by the relation
\begin{align}
p_h = \beta^*(\hat{p}_h - \bar{w}_h) \mathcal{I}_h^k(\Psi G) + v_h, \label{eq:reconstruct}
\end{align}
where $\mathcal{I}_h^k$ denotes the interpolation operator onto the Lagrange elements of order $k$. Finally, the numerical error associated with $v_h$ can be computed as 
\begin{align}
p_e = I_h^{k=1}(p_a) - p_h,
\end{align}
where $p_a$ is interpolated onto the Lagrange elements with the same order as the solution $v_h$.

The results of applying the SRB-FE method to solve this problem are plotted in Figure \ref{fig:exp2-pressure} for $h = \sfrac{1}{8}$. The errors and convergence rates are given for different mesh sizes in Table \ref{tab:exp2}. As is evident from this table, the SRB-FE approximation of $p_h$ and $\hat{p}_h$ both converge with optimal order. I.e., we find that
\begin{align}
\norm{p_e}{L^2(\Omega)} &\leq C h^2 \norm{I_h^{k=1}(p_a)}{H^1(\Omega)},  \\
\norm{p_e}{H^1(\Omega)} &\leq C h^1 \norm{I_h^{k=1}(p_a)}{H^2(\Omega)},  \\
\norm{\hat{p}_e}{L^2(\Lambda)} &\leq C h^2 \norm{\hat{p}}{H^1(\Lambda)}, \\
\norm{\hat{p}_e}{H^1(\Lambda)} &\leq C h^1 \norm{\hat{p}}{H^2(\Lambda)}.
\end{align}

\section{Conclusion} 
 
In this work, we have developed a singularity removal method for the coupled 1D-3D flow model. This type of model can be used to model the interaction of wells with a reservoir. The well is endowed with its own 1D flow equation, and modelled as a 1D line source in the reservoir domain. This line source introduces a logarithmic type singularity in the reservoir solution that negatively affects the approximation properties of the problem. We provide here a method for identifying and removing this singularity from the governing equations. The result is a reformulated coupled 1D-3D flow model in which all variables are smooth. 

As the reformulated model is posed in terms of smooth variables, it has the advantage that it can be approximated using any standard numerical method. In this work, we have shown that the singularity removal restores optimal convergence rates for the Galerkin FE method. Moreover, it makes the approximation stable with respect to the ratio $R/h$ between well radius and mesh size. 

A natural development of this work consists of extending the singularity removal method to apply to (i) different control modes for the wells, (ii) tensor-valued permeability and (iii) a mixed formulation of the flow, where both pressure and flux are approximated. We believe these extensions would be particularly valuable in the context of subsurface flow applications, as it would allow one to capture the interaction between well and reservoir using coarse grids. The extension to different control modes for the wells, i.e., rate controlled or pressure controlled wells, is straightforward; it can be achieved by altering the boundary conditions for the well flow equations. As the singularity subtraction is performed at the continuous level, it is likewise straightforward to adapt the method to different discretization methods \cite{ecmor}. The extension to tensor-valued permeability is more challenging, and will be treated in future work. 

\vspace{5em}

\begin{acknowledgements}
This work was partially supported by the Research Council of Norway, project number 250223. The participation of the first author at ECMOR XVI 2018 was funded by the Academia Agreement between the University of Bergen and Equinor. 

The authors thank B. Wohlmuth for her help in formulating the original splitting approach. The authors also thank P. Zunino for the interesting discussions on coupled 1D-3D flow problems. Finally, the authors would like to thank M. Kuchta, K. E. Holter and C. Daversin-Catty for developing and sharing with us code for the implementation of mixed-dimensional models in FEniCS.
\end{acknowledgements}

\bibliographystyle{spmpsci}      

\bibliography{references}

\begin{thebibliography}{10}
\providecommand{\url}[1]{{#1}}
\providecommand{\urlprefix}{URL }
\expandafter\ifx\csname urlstyle\endcsname\relax
  \providecommand{\doi}[1]{DOI~\discretionary{}{}{}#1}\else
  \providecommand{\doi}{DOI~\discretionary{}{}{}\begingroup
  \urlstyle{rm}\Url}\fi

\bibitem{Aavatsmark2016-0}
Aavatsmark, I.: Equivalent well-cell radius for hexagonal k-orthogonal grids in
  numerical reservoir simulation.
\newblock Applied Mathematics Letters \textbf{61}, 122 -- 128 (2016)

\bibitem{Aavatsmark2016}
Aavatsmark, I.: Interpretation of well-cell pressures on hexagonal grids in
  numerical reservoir simulation.
\newblock Computational Geosciences \textbf{20}(5), 1029--1042 (2016)

\bibitem{Aavatsmark2016-2}
Aavatsmark, I.: Interpretation of well-cell pressures on stretched hexagonal
  grids in numerical reservoir simulation.
\newblock Computational Geosciences \textbf{20}(5), 1043--1060 (2016)

\bibitem{aavatsmark2003index}
Aavatsmark, I., Klausen, R.A.: Well index in reservoir simulation for slanted
  and slightly curved wells in 3d grids.
\newblock SPE Journal pp. 41--48 (2003).
\newblock \doi{10.2118/75275-PA}

\bibitem{alkhoury2005}
Al$-$Khoury, R., Bonnier, P.G., Brinkgreve, R.B.J.: Efficient finite element
  formulation for geothermal heating systems. part i: steady state.
\newblock International Journal for Numerical Methods in Engineering
  \textbf{63}(7), 988--1013 (2005)

\bibitem{apel2011}
Apel, T., Benedix, O., Sirch, D., Vexler, B.: A priori mesh grading for an
  elliptic problem with dirac right-hand side.
\newblock SIAM Journal on Numerical Analysis \textbf{49}(3), 992--1005 (2011)

\bibitem{babu}
Babu, D.K., Odeh, A.S., Al-Khalifa, A.J., McCann, R.C.: The relation between
  wellblock and wellbore pressures in numerical simulation of horizontal wells
  (1991)

\bibitem{babuska1972}
Babu{\v{s}}ka, I., Rosenzweig, M.B.: A finite element scheme for domains with
  corners.
\newblock Numerische Mathematik \textbf{20}(1), 1--21 (1972)

\bibitem{splineinterpolation}
Broomhead, D.S., Lowe, D.: {Multivariable Functional Interpolation and Adaptive
  Networks}.
\newblock Complex Systems 2 pp. 321--355 (1988)

\bibitem{Cattaneo2014}
Cattaneo, L., Zunino, P.: A computational model of drug delivery through
  microcirculation to compare different tumor treatments.
\newblock International Journal for Numerical Methods in Biomedical Engineering
  \textbf{30}(11), 1347--1371 (2014)

\bibitem{Cerroni2019}
Cerroni, D., Laurino, F., Zunino, P.: Mathematical analysis, finite element
  approximation and numerical solvers for the interaction of 3d reservoirs with
  1d wells.
\newblock GEM - International Journal on Geomathematics \textbf{10}(1), 4
  (2019)

\bibitem{Chen2009}
Chen, Z., Zhang, Y.: Well flow models for various numerical methods.
\newblock International Journal of Numerical Analysis and Modeling \textbf{3},
  375--388 (2006)

\bibitem{dangelo2012}
D'Angelo, C.: Finite element approximation of elliptic problems with dirac
  measure terms in weighted spaces: Applications to one- and three-dimensional
  coupled problems.
\newblock SIAM Journal on Numerical Analysis \textbf{50}(1), 194--215 (2012)

\bibitem{dangelo2008}
D'Angelo, C., Quarteroni, A.: On the coupling of 1d and 3d diffusion-reaction
  equations: Application to tissue perfusion problems.
\newblock Mathematical Models and Methods in Applied Sciences \textbf{18}(08),
  1481--1504 (2008)

\bibitem{cecile-git}
Daversin-Catty, C.:
  \url{https://hub.docker.com/r/ceciledc/fenics_mixed_dimensional/}

\bibitem{Ding2001}
Ding, Y., Jeannin, L.: A new methodology for singularity modelling in flow
  simulations in reservoir engineering.
\newblock Computational Geosciences \textbf{5}(2), 93--119 (2001)

\bibitem{ewing-book}
Ewing, R.E. (ed.): The mathematics of reservoir simulation, \emph{Frontiers in
  Applied Mathematics}, vol.~1.
\newblock Society for Industrial and Applied Mathematics (SIAM), Philadelphia,
  PA (1983)

\bibitem{Ewing1999}
Ewing, R.E., Lazarov, R.D., Lyons, S.L., Papavassiliou, D.V., Pasciak, J., Qin,
  G.: Numerical well model for non-darcy flow through isotropic porous media.
\newblock Computational Geosciences \textbf{3}(3), 185--204 (1999)

\bibitem{ecmor}
Gjerde, I.G., Kumar, K., Nordbotten, J.M.: Well modelling by means of coupled
  1d-3d flow models (2018).
\newblock \doi{10.3997/2214-4609.201802117}

\bibitem{Gjerde2018}
Gjerde, I.G., Kumar, K., Nordbotten, J.M., Wohlmuth, B.: {Splitting method for
  elliptic equations with line sources}.
\newblock ArXiv e-prints p. arXiv:1810.12979 (2018)

\bibitem{Grinberg2011}
Grinberg, L., Cheever, E., Anor, T., Madsen, J.R., Karniadakis, G.E.: Modeling
  blood flow circulation in intracranial arterial networks: A comparative 3d/1d
  simulation study.
\newblock Annals of Biomedical Engineering \textbf{39}(1), 297--309 (2011)

\bibitem{hales}
Hales, H.B.: An improved method for simulating reservoir pressures through the
  incorporation of analytical well functions  (1977).
\newblock \doi{10.2118/39065-MS}

\bibitem{Holter2018}
{Holter}, K.E., {Kuchta}, M., {Mardal}, K.A.: {Sub-voxel perfusion modeling in
  terms of coupled 3d-1d problem}.
\newblock ArXiv e-prints  (2018)

\bibitem{king}
King, M.J., Mansfield, M.: Flow simulation of geologic models  (1997).
\newblock \doi{doi:10.2118/39065-MS}

\bibitem{koppl2016}
K{\"o}ppl, T., Vidotto, E., Wohlmuth, B.I., Zunino, P.: Mathematical modelling,
  analysis and numerical approximation of second order elliptic problems with
  inclusions.
\newblock Numerical Mathematics and Advanced Applications ENUMATH 2015  (2017)

\bibitem{miro-git}
Kuchta, M.: \url{https://github.com/MiroK}, repo: fenics$_{}$ii, branch: master

\bibitem{miro2016-3D1D}
{Kuchta}, M., {Mardal}, K.A., {Mortensen}, M.: {Preconditioning trace coupled
  3$d$-1$d$ systems using fractional Laplacian}.
\newblock ArXiv e-prints  (2016)

\bibitem{koppl2015}
Köppl, T., Vidotto, E., Wohlmuth, B.: A local error estimate for the poisson
  equation with a line source term.
\newblock Numerical Mathematics and Advanced Applications ENUMATH pp. 421--429
  (2015)

\bibitem{llau2016}
Llau, A., Jason, L., Dufour, F., Baroth, J.: Finite element modelling of 1d
  steel components in reinforced and prestressed concrete structures.
\newblock Engineering Structures \textbf{127}(Supplement C), 769 -- 783 (2016)

\bibitem{LoggMardalEtAl2012a}
Logg, A., Mardal, K.A., Wells, G.N.: Automated Solution of Differential
  Equations by the Finite Element Method.
\newblock Springer (2012).
\newblock \doi{10.2118/75275-PA}

\bibitem{nabil2016}
Nabil, M., Zunino, P.: A computational study of cancer hyperthermia based on
  vascular magnetic nanoconstructs.
\newblock Royal Society Open Science \textbf{3}(9) (2016)

\bibitem{Nordbotten2009}
Nordbotten, J.M., Kavetski, D., Celia, M.A., Bachu, S.: Model for co2 leakage
  including multiple geological layers and multiple leaky wells.
\newblock Environmental Science \& Technology \textbf{43}(3), 743--749 (2009).
\newblock PMID: 19245011

\bibitem{peaceman1978}
Peaceman, D.W.: Interpretation of well-block pressures in numerical reservoir
  simulation.
\newblock Society of Petroleum Engineers Journal \textbf{18}(03), 183--194
  (1978)

\bibitem{Peaceman-rectangle}
Peaceman, D.W.: Interpretation of well-block pressures in numerical reservoir
  simulation with nonsquare grid blocks and anisotropic permeability.
\newblock Society of Petroleum Engineers Journal \textbf{23}(3) (1983)

\bibitem{zunino2018}
Possenti, L., Casagrande, G., Gregorio, S.D., Zunino, P., Constantino, M.:
  {Numerical simulations of the microvascular fluid balance with a non-linear
  model of the lymphatic system}.
\newblock MOX-Report No. 35  (2018)

\bibitem{Reichold2009}
Reichold, J., Stampanoni, M., Keller, A.L., Buck, A., Jenny, P., Weber, B.:
  Vascular graph model to simulate the cerebral blood flow in realistic
  vascular networks.
\newblock Journal of Cerebral Blood Flow \& Metabolism \textbf{29}(8),
  1429--1443 (2009)

\bibitem{gfem}
Strouboulis, T., Babu\v{s}ka, I., Copps, K.: The design and analysis of the
  generalized finite element method.
\newblock Comput. Methods Appl. Mech. Engrg. \textbf{181}(1-3), 43--69 (2000)

\bibitem{weiss2017}
Weiss, C.J.: Finite-element analysis for model parameters distributed on a
  hierarchy of geometric simplices.
\newblock GEOPHYSICS \textbf{82}(4), E155--E167 (2017)

\bibitem{Wolfsteiner2003}
Wolfsteiner, C., Durlofsky, L.J., Aziz, K.: Calculation of well index for
  nonconventional wells on arbitrary grids.
\newblock Computational Geosciences \textbf{7}(1), 61--82 (2003)

\bibitem{hamdi2006}
Wolfsteiner, C., Lee, S., Tchelepi, H.: Well modeling in the multiscale finite
  volume method for subsurface flow simulation.
\newblock Multiscale Modeling \& Simulation \textbf{5}(3), 900--917 (2006)

\end{thebibliography}

\end{document}